\documentclass[11pt]{article}
\usepackage{graphicx}
\usepackage{amssymb}
\usepackage{amsmath}
\usepackage{subfigure}
\usepackage{slashbox}
\usepackage{multirow}

\renewcommand{\u}{{\bf u}}
\newcommand{\x}{{\bf x}}
\newcommand{\f}{{\bf f}}
\newcommand{\F}{{\bf F}}
\renewcommand{\S}{{\bf S}}
\newtheorem{theorem}{Theorem}

\title{Termination criteria for inexact fixed point schemes}
\author{Philipp Birken$^{\mbox{\tiny\rm 1}}$}
\begin{document}
\maketitle
\baselineskip=0.9
\normalbaselineskip
\vspace{-3pt}
\begin{center}{\footnotesize\em $^{\mbox{\tiny\rm 1}}$
Numerical Analysis, Centre for the Mathematical Sciences, Lund University, Box 118, 22100 Lund, Sweden\\
Institute of Mathematics, University of Kassel, Heinrich-Plett-Str. 40, 34132 Kassel, Germany\\
Department of Mathematics/Computer Science, University of Osnabr\"uck, Albrechtstr. 28a, 49076 Osnabr\"uck, Germany\\
email: philipp.birken\symbol{'100}na.lu..de}
\end{center}

\begin{abstract}
We analyze inexact fixed point iterations where the generating
function contains an inexact solve of an equation system to answer the question of how tolerances for the inner solves influence the iteration error of the outer fixed point iteration. Important 
applications are the Picard iteration and partitioned fluid structure
interaction. We prove that the iteration converges irrespective of how
accurate the inner systems are solved, provided that a specific
relative termination criterion is employed, whereas standard relative and
absolute criteria do not have this property. For the analysis, the iteration is modelled as a perturbed fixed point
iteration and existing analysis is extended to the nested case
$\x={\bf F}({\bf S}(\x))$. 
\end{abstract}

{\it Keywords: Fixed point iteration, Picard iteration, Transmission
  Problem, Dirichlet-Neumann iteration, Termination criteria}

\section{Introduction}

The general problem that this article is about is the following: Consider a nonlinear equation system and an outer iteration method to solve it that consists of solving a subproblem at each step using a second, inner iteration method. Now we want to answer the following question: How can we efficiently control the iteration error of the outer iteration method? Or otherwise put: How accurate do we need to solve the inner systems to obtain a certain iteration error for the outer nonlinear equation? 

For the case of the outer iteration being Newton's method, this problem has been successfully solved. The inner problem is a linear system and the concept of an inexact Newton's method was introduced in \cite{deeist:82}. There, at each Newton step the inner iteration is terminated when a relative tolerance criterion in the linear residual is satisfied. Based on this, it is possible to give conditions on the sequence of relative tolerances to obtain linear, superlinear or quadratic convergence of the inexact Newton's method. Essentially, the sequence of tolerances has to converge to zero fast enough as the Newton scheme progresses and then quadratic convergence is obtained. Following up, a strategy that has this property and leads to a very efficient scheme was suggested in \cite{eiswal:96}. There, the point is that the initial systems are solved quite coarsely and these schemes are part of widely used software packages, for example of SUNDIALS \cite{hbglss:05}.

Note that with this knowledge, the choice of iterative solver for the inner iteration obtains a better basis: If most of the systems are solved very coarsely and thus very few iterations are needed, it is more important that the method is cheap per iteration than how fast we can reach machine accuracy. In this setting, when looking at unsymmetric linear systems and Krylov subspace methods, GMRES \cite{saadsch:86} beats BiCGSTAB \cite{vorst:92}, since it needs only one matrix vector product per iteration instead of two. 

Now when looking at fixed point iterations, an iteration typically consists of evaluating a function and not of solving a system. However, two prominent and important examples where this happens are the Picard iteration and fluid-structure interaction and thus we call these inexact fixed point schemes. Surprisingly, the problem framed in the first paragraph has not been analyzed for these. 

For the Picard iteration, which is a common tool in the context of the incompressible Navier-Stokes equations, the evaluation of the right hand side corresponds to solving a linear system. Strategies for choosing a termination criterion for the inner iteration are empirically discussed for example in \cite{elsiwa:05,turek:99}. 

In Fluid-Structure interaction, a standard approach are partitioned coupling schemes, where existing solvers for the subproblems are reused \cite{farhat:04}. Commonly in the form of a Dirichlet-Neumann iteration, this consists of subsequently solving the fluid and the structure problem with appropriate boundary conditions and reasonable tolerances. It is common to formulate the coupling condition at the interface in the form of a fixed point equation. Recently, it was suggested to use a time adaptive implicit time integration scheme for fluid structure interaction \cite{biquhm:11}, where the time step is chosen based on an error tolerance. As is common in this setting, the tolerances for the solvers for the appearing nonlinear equation systems are chosen such that the iteration error does not interfere with the error from the time integration scheme \cite{birkenhabil}, but nevertheless as large as possible to avoid unnecessary computations. Thus, it is imperative to be able to control the iteration error. 

To solve our problem, we proceed in the following way. First, we will review well known results on perturbations of fixed point schemes \cite{ortrhe:00}. The general idea is then to quantify the iteration errors based on the termination criterion of the inner iteration such that the existing perturbation results can be applied. For the Dirichlet-Neumann iteration, we first have to extend the perturbation results to the case of a nested fixed point equation of the form $\x={\bf F}({\bf S}(\x))$. 

As it turns out, the type of termination criterion chosen in the inner solver is crucial for the answer to our problem in that when using a nonstandard relative criterion, we obtain convergence of the fixed point iteration to the exact solution independently of how accurate we solve the inner systems, leading to an efficient way of controlling the outer iteration error. 

On the other hand, no general statement can be made on a standard relative or absolute termination criterion, but the analysis suggests that these do not have favorable properties. All of this is confirmed by numerical results, which show that the latter criteria cause convergence to a solution that is farther away from the exact one the less accurate we solve the linear systems or otherwise put, the less error we want in the fixed point equation, the more accurate we have to solve the linear systems.

\section{Inexact Fixed Point Methods}

\subsection{Direct perturbation}

Consider the fixed point equation
\begin{equation}\label{fixedpointeq}
\x=\f(\x)
\end{equation}
with $\x \in \Omega \subset \mathbb{R}^n$, $\Omega$ closed and where
we assume that $\f:\Omega \rightarrow \Omega$ is Lipschitz
continuous with Lipschitz constant $L<1$. This implies by the Banach
fixed point theorem that
\eqref{fixedpointeq} has a unique solution $\x^*$. 

Furthermore, we consider the perturbed fixed point iteration
\begin{equation}\label{perturbedfixedpointiter}
\x^{k+1}=\f(\x^k) + \epsilon,
\end{equation}
where $\epsilon$ is a perturbation that could originate from an iterative solver and for simplicities sake we denote the norm of $\epsilon$ by $\epsilon$ as well. We furthermore assume for simplicities sake that $\f + \epsilon$ is also a self-map on $\Omega$. Thus this iteration obtains a solution $\x_{\epsilon}$ of the perturbed fixed point equation 
\begin{equation}\label{perturbedfixedpointeq}
\x=\f(\x) + \epsilon.
\end{equation}
The question is now: How far is the solution $\x_{\epsilon}$ of that equation away from $\x^*$? The answer is giving by the following theorem, see for example \cite{ortrhe:00}. 
\begin{theorem}\label{ortrhetheorem}
For the solutions $\x_{\epsilon}$ and $\x^*$ of problems \eqref{perturbedfixedpointeq} and \eqref{fixedpointeq} we have:
\begin{equation}\label{xepsminusxstar}
\|\x_{\epsilon}-\x^*\| \leq \epsilon \frac{1}{1-L}.
\end{equation}
\end{theorem}

This means that the error is of the order $\epsilon$ as is to be expected, but interestingly, it becomes larger, the closer the Lipschitz constant of $\f$ is to one or otherwise put, the less contractive the function is. This implies that in these cases, the error will be much larger than $\epsilon$ and thus a much smaller
tolerance would have to be supplied to acchieve the desired error. Note that in practice, we typically do not know the Lipschitz constant and that in the nonlinear case, it depends on the definition of the domain $\Omega$. Thus the important Lipschitz constant is the local one in the solution.  

If we instead consider a sequence of perturbations $\epsilon_k$, respectively a nonconstant perturbation, and thus the iteration
\begin{equation}\label{perturbedfixedpoint}
\x^{k+1}=\f(\x^k)+\epsilon_k,
\end{equation} 
the first question is when this sequence converges to $\x^*$. The answer is given by the next theorem, also from \cite{ortrhe:00}:
\begin{theorem}\label{perturbedtheorem}
The iteration \eqref{perturbedfixedpoint} converges to the solution of the unperturbed problem \eqref{fixedpointeq} if and only if $\lim_{k \rightarrow \infty}\epsilon_k=0$. 
\end{theorem}
A specific case is
\begin{equation}\label{adaptivestrat}
\epsilon_k=cL^k
\end{equation}
with $c>0$, which we call the adaptive strategy.

\subsection{Application: Picard iteration}

As an application of the above theorems, we now analyze the convergence of the Picard iteration. This is often employed in the context of the incompressible Navier-Stokes equation and corresponds to a fixed point
iteration for the equation
\[{\bf x}={\bf A}^{-1}({\bf x}){\bf b},\]
where ${\bf A}({\bf x})\in \mathbb{R}^{n\times n}$ is an approximation of a
Jacobian in ${\bf x}$ \cite{elsiwa:05}. Thus, the fixed point iteration
\begin{equation}\label{picardfixed}
{\bf x}^{k+1}={\bf A}^{-1}({\bf x}^k){\bf b}
\end{equation}
is implemented by solving
\begin{equation}\label{picardlineq}
{\bf A}({\bf x}^k){\bf x}^{k+1}={\bf b}
\end{equation}
for $\x^{k+1}$ up to a certain tolerance using an iterative scheme. The scheme
\eqref{picardfixed} can be analyzed either as a fixed point scheme,
which results in linear convergence provided that the Lipschitz
constant $L$ of ${\bf A}^{-1}(\x){\bf b}$ can be bounded from above away from one or as a method of
Newton type where ${\bf A}(\x)$ is an approximation of
the exact Jacobian and we have linear convergence as long as this
approximation is good enough. 

When solving \eqref{picardlineq}, either the relative termination criterion
\begin{equation}\label{linrelterm}
\|{\bf A}({\bf x}^k){\bf x}^{k+1}-{\bf b}\| \leq \tau_r \|{\bf A}({\bf x}^k){\bf x}^k-{\bf b}\|,
\end{equation}
the relative termination criterion
\begin{equation}\label{linreltermb}
\|{\bf A}({\bf x}^k){\bf x}^{k+1}-{\bf b}\| \leq \tau_r \|{\bf b}\|
\end{equation}
or the absolute criterion
\begin{equation}\label{linabsterm}
\|{\bf A}({\bf x}^k){\bf x}^{k+1}-{\bf b}\| \leq \tau_a
\end{equation}
are used, where $\tau_r$ and $\tau_a$ are relative and adaptive tolerances. 

To analyze the consequences of choosing one of these using theorems 1 and 2, we need to quantify the perturbation in the form
\eqref{perturbedfixedpointiter}. Thus, we define ${\bf f}(\x^k)={\bf
  A}^{-1}(\x^k){\bf b}$ to obtain
\[\x^{k+1}={\bf A}^{-1}(\x^k){\bf b} + \epsilon_k\]
and we can write 
\begin{equation}\label{xepsilon}
\epsilon_k={\bf A}^{-1}(\x^k)({\bf A}({\bf x}^k){\bf x}^{k+1}-{\bf b}).
\end{equation}

In the case of the relative termination criterion \eqref{linrelterm},
we can estimate the norm of the right hand side in \eqref{xepsilon} by
\[\|\epsilon_k\|\leq \|{\bf A}^{-1}(\x^k)\|\tau_r\|({\bf A}({\bf
  x}^k){\bf x}^k-{\bf b})\|.\]
We furthermore have
\[{\bf A}({\bf x}^k){\bf x}^k-{\bf b}={\bf A}({\bf x}^k)({\bf
  x}^k-{\bf A}^{-1}({\bf x}^k){\bf b})={\bf A}({\bf x}^k)({\bf
  f}(\x^{k-1})-{\bf f}(\x^k)).\]
Thus, 
\[\|({\bf A}({\bf x}^k){\bf x}^k-{\bf b})\| \leq \|({\bf A}({\bf
  x}^k)\|L\|\x^{k-1}-\x^k\| \leq \|({\bf A}({\bf
  x}^k)\|L^k\|\x^1-\x^0\|.\]
All in all, we obtain with the condition number $\kappa({\bf A})$
\begin{equation}
\|\epsilon_k\|\leq \tau_r \kappa({\bf
    A}(\x^k)) L^k \|\x^1-\x^0\|.
\end{equation}
 With the additional and reasonable assumption that $\kappa({\bf
  A}(\x))$ is bounded, this is a perturbation of the form \eqref{adaptivestrat} and $\epsilon_k$ converges to zero independent
of the choice of $\tau_r$! Thus by theorem \ref{perturbedtheorem} this iteration converges to the exact solution independently of how accurate we solve the linear equation systems. We
now formulate this as a theorem. 
\begin{theorem}\label{picardtheorem}
Let ${\bf b}\in \mathbb{R}^n$ and the function ${\bf A}(\x)$ be given
that maps the closed set $\Omega \subset \mathbb{R}^n$ onto quadratic
regular matrices. Assume that the function ${\bf A}^{-1}(\x){\bf b}:\Omega\rightarrow
\Omega$ is Lipschitz continuous with Lipschitz constant $L<1$ and correspondingly has a unique fixpoint $\x^*$. Furthermore assume that $\kappa({\bf A}(\x))$ is bounded on $\Omega$ and that the inexact fixedpoint iteration defined by
\eqref{linrelterm} converges to a limit $\x_{\epsilon}$. 
Then $\x_{\epsilon}=\x^*$, independent of the choice of $\tau_r$. 
\end{theorem}

In case of the relative criterion \eqref{linreltermb}, the estimate
\[\|\epsilon_k\|\leq \|{\bf A}^{-1}(\x^k)\|\tau_r \|{\bf b}\|\]
for the norm of the left hand side in \eqref{xepsilon} holds which is bounded away from zero provided that ${\bf A}^{-1}(\x)$
is. Thus, it is not clear if this iteration satisfies theorem
\ref{perturbedtheorem}, but in the general case, the iteration
will not converge to $\x^*$. Similarly if we use the absolute termination criterion \eqref{linabsterm}, we
obtain
\[\|\epsilon_k\|\leq \|{\bf A}^{-1}(\x^k)\|\tau_a\]
which is also bounded away from zero if ${\bf A}^{-1}(\x)$ is. 

Numerical results that confirm theorem \ref{picardtheorem} and
demonstrate that the other two iterations behave like being of the form
\eqref{ortrhetheorem} can be found in section
\ref{picardresultssection}. We would like to point out that the
criterion \eqref{linrelterm} is sometimes suggested in the literature on the Picard iteration, e.g. \cite{johmat:04, elsiwa:05}, but that an absolute termination criterion is suggested in \cite{turek:99}. There it is suggested to just ``gain one digit'', meaning to use a tolerance of 0.1.

\subsection{Perturbed nested fixed point iteration}

Now consider two functions $\F:\Omega_1 \rightarrow \Omega_2$ and $\S:\Omega_2\rightarrow \Omega_1$ with $\Omega_1, \Omega_2 \subset \mathbb{R}^n$ closed and the fixed point equation
\begin{equation}
\x=\S(\F(\x))
\end{equation}
again with solution $\x^*$. We
now consider an iteration where both the evaluation of $\F$
and of $\S$ are perturbed, namely $\S$ is perturbed by $\delta_k$ and
$\F$ by $\epsilon_k$:
\begin{equation}\label{nestedperturbedeq}
\x^{k+1}=\S(\F(\x^k)+\epsilon_k)+\delta_k.
\end{equation}
Again, assume that this iteration is well defined and that this sequence has the limit $\x_{\epsilon}$. Then, we obtain the following theorem. 
\begin{theorem}\label{fsitheorem}
Let $\F$ and $\S$ be Lipschitz
continuous with Lipschitz constants $L_F$ and $L_S$,
respectively. Assume that $L_FL_S<1$. Then we have, if $\epsilon_k=\delta_k=\epsilon$ for
all $k$, that
\begin{equation}\label{fsiestimateepseqdelta}
\|\x_{\epsilon}-\x^*\|\leq \epsilon
\frac{1+L_S}{1-L_SL_F}.
\end{equation}
In the case $\epsilon_k=\epsilon$ and $\delta_k=\delta$, we obtain
\begin{equation}\label{fsiestimateepsdelta}
\|\x_{\epsilon}-\x^*\|\leq \frac{\epsilon L_S+\delta}{1-L_SL_F}.
\end{equation}
Finally, $\x_{\epsilon}=\x^*$ if and only if both $\delta_k$ and $\epsilon_k$
converge to zero. 
\end{theorem}

Proof: The proof is technically identical to the one of theorem
\ref{ortrhetheorem}. We have due to the Lipschitz continuity
\begin{eqnarray*}
\|\x^{k+1}-\x^*\|=\|\S(\F(\x^k)+\epsilon_k)+\delta_k-\x^*\|=\|\S(\F(\x^k)+\epsilon_k)+\delta_k
- \S(\F(\x^*))\|\\
\leq L_S\|\F(\x^k)-\F(\x^*)+\epsilon_k\|+\delta_k \leq
L_SL_F\|\x^k-\x^*\| + L_S\epsilon_k + \delta_k \\
\leq
(L_SL_F)^2\|\x^{k-1}-\x^*\| + L_S^2L_F\epsilon_{k-1} + L_SL_F\delta_{k-1} +
L_S\epsilon_k + \delta_k \\
\leq (L_SL_F)^{k+1}\|\x^0-\x^*\| + \left(
  \sum_{j=0}^kL_S^{j+1}L_F^j \epsilon_{k-j}\right ) +  \left(
  \sum_{j=0}^kL_S^jL_F^j \delta_{k-j}\right )
\end{eqnarray*}
and thus in the limit $\x^{k+1} \rightarrow \x^*$,
\begin{equation}\label{nestedgeneral}
\|\x_{\epsilon}-\x^*\|\leq L_S \lim_{k\rightarrow \infty}
\sum_{j=0}^k(L_SL_F)^j \epsilon_{k-j} + \lim_{k\rightarrow \infty}
\sum_{j=0}^k(L_SL_F)^j \delta_{k-j}
\end{equation}

For a constant perturbation overall, e.g. $\epsilon_k=\delta_k=\epsilon$ for
all $k$, we obtain in the limit
\[\|\x_{\epsilon}-\x^*\|\leq \epsilon (1+L_S) \lim_{k\rightarrow \infty}
\sum_{j=0}^k(L_SL_F)^j
=\epsilon
\frac{1+L_S}{1-L_SL_F},\]
which proves the inequality \eqref{fsiestimateepseqdelta}. If we have constant but separate perturbations $\epsilon$ and $\delta$
of $\S$ and $\F$, we obtain \eqref{fsiestimateepsdelta} from
\[\|\x_{\epsilon}-\x^*\|\leq \epsilon L_S \lim_{k\rightarrow \infty}
\sum_{j=0}^k(L_SL_F)^j + \delta \lim_{k\rightarrow \infty}
\sum_{j=0}^k(L_SL_F)^j
=
\frac{\epsilon L_S+\delta}{1-L_SL_F}.\]

In the general case, due to positivity, the right hand side of
\eqref{nestedgeneral} is zero if and only if both $\epsilon_k$ and
$\delta_k$ are such that for $\phi_k=\epsilon_k$ or $\phi_k=\delta_k$,
\[\lim_{k\rightarrow \infty} \sum_{j=0}^k(L_SL_F)^j \phi_{k-j}=0.\]
By an identical proof to theorem 2, this is the case if and only if
both $\epsilon_k$ and $\delta_k$ converge to zero. 

Note that this implies that the sequence $\epsilon_k$ perturbing the
inner function $\F$ is less important by a factor of the Lipschitz
constant $L_S$ of the outer function. Thus, a possible strategy is to
define
\begin{equation}\label{adaptivestrat2}
\epsilon_k=\delta_k/L_S,
\end{equation}
meaning that we solve the fluid part less accurate by a factor of
$L_S$. Unfortunately, $L_S$ has to be known for this.

\subsection{Application: Dirichlet-Neumann coupling for Transmission problem}

As an application of the theory from section 2.3, we consider a problem that is a basic building block in fluid structure
interaction, namely the transmission problem, where the Laplace
equation with right hand side $f(x,y)$ on a domain $\Omega$ is cut into two domains $\Omega=\Omega_1 \cup \Omega_2$ using transmission
conditions at the interface $\Gamma=\Omega_1\cap
\Omega_2$:
\begin{eqnarray}\label{transmissioneq}
\Delta u_i(x,y)&=f(x,y), \, (x,y) \in \Omega_i \subset \mathbb{R}^2, \, i=1,2 \nonumber \\
u_i(x,y)&=0, \, (x,y)\in \partial \Omega_i \ \partial \Omega_1\cap
\Omega_2\\
u_1(x,y)&=u_2(x,y), \, (x,y) \in \Gamma \nonumber\\
\partial u_1(x,y)\cdot {\bf n}&=\partial u_2(x,y)\cdot {\bf n}, \,
(x,y) \in \Gamma \nonumber
\end{eqnarray}

We now employ a standard Dirichlet-Neumann
iteration to solve it. Using any linear discretization, this corresponds to alternately solving the problems
\begin{equation}\label{dirichletproblem}
{\bf A u}_1^{k+1}={\bf b}_1({\bf u}_2^k)
\end{equation}
and 
\begin{equation}\label{neumannproblem}
{\bf B u}_2^{k+1}={\bf b}_2({\bf u}_1^{k+1})
\end{equation}
were problem \eqref{dirichletproblem} corresponds to a discretization
of the transmission
problem \eqref{transmissioneq} on $\Omega_1$ only with Dirichlet data
on $\Gamma$ given by ${\bf u}_2^k$ on the
coupling interface and problem \eqref{neumannproblem} corresponds to a
discretization of \eqref{transmissioneq} on $\Omega_2$ only with
Neumann data on $\Gamma$ given by the discrete normal derivative of
${\bf u}_1$ on $\Gamma$. It can be shown that
convergence of the approximate solutions on the whole domain is
equivalent to the convergence of the solution on the interface only
\cite{quaval:99}. 

By considering \eqref{dirichletproblem}-\eqref{neumannproblem} as one iteration, we
obtain a fixed point formulation
\[\u_{\Gamma}=\S(\F({\bf u}_{\Gamma}))\]
where $\u_{\Gamma}$ is ${\bf u}_2$ on the interface, $\F={\bf D}_{{\bf n}_{\Gamma}}{\bf
  A}^{-1}{\bf b}_1({\bf u}_{\Gamma})$ and $\S={\bf P}_{\Gamma}{\bf
  B}^{-1}{\bf b}_2({\bf u}_1)$. Hereby ${\bf D}_{{\bf n}_{\Gamma}}$ is
the matrix that computes the discrete normal derivatives in $\Omega_1$ on $\Gamma$
and ${\bf P}_{\Gamma}$ is the discrete trace operator with respect to $\Gamma$. Otherwise put, ${\bf P}_{\Gamma}$ is the projection of the space that ${\bf u}_2$
is in onto the space of discrete unknowns on $\Gamma$. 

In practice, the linear equation systems are
solved iteratively, typically using the conjugate gradient method (CG) up to a relative tolerance of $\tau$. Thus,
we obtain a perturbed nested fixed point iteration of the form
\eqref{nestedperturbedeq} and the question is now again if we can
quantify this perturbation. We have
\begin{equation}\label{epsilonbound}
{\bf u}_{1_{\epsilon}}^{k+1}={\bf A}^{-1}{\bf b}_1({\bf u}_{\Gamma}^k)
+ \epsilon_k
\end{equation}
and
\begin{equation}\label{deltabound}
{\bf u}_{2_{\epsilon}}^{k+1}={\bf B}^{-1}{\bf b}_2({\bf u}_{1_{\epsilon}}^{k+1})
+\delta_k.
\end{equation}

For the iteration \eqref{epsilonbound} we obtain
\[\|\epsilon_k\|=\|{\bf u}_{1_{\epsilon}}^{k+1}-{\bf
  A}^{-1}{\bf b}_1({\bf u}_{\Gamma}^k)\| \leq \|{\bf A}^{-1}\|\|{\bf
  A}{\bf u}_{1_{\epsilon}}^{k+1}-{\bf b}_1({\bf u}_{\Gamma}^k)\|.\]
Again, the second factor is what is tested in the termination
criterion of CG. In the case of the relative criterion
\eqref{linrelterm}, here stated as
\[\|{\bf A}{\bf u}_{1_{\epsilon}}^{k+1}-{\bf b}_1({\bf
  u}_{\Gamma}^k)\|\leq \tau_r \|{\bf A}{\bf u}_1^k-{\bf b}_1({\bf u}_{\Gamma}^k)\|,\]
we obtain
\[\|\epsilon_k\| \leq \|{\bf A}^{-1}\|\tau_r \|{\bf
  A}{\bf u}_{1_{\epsilon}}^k-{\bf b}_1({\bf u}_{\Gamma}^k)\|\leq
\kappa({\bf A})\tau_r \|{\bf u}_1^k-{\bf A}^{-1}{\bf b}_1({\bf
  u}_{\Gamma}^k)\|\]
\[= \kappa({\bf A}) \tau_r \|{\bf u}_1^k-{\bf u}_1^{k+1}\|.\]

Now the point is that since the whole iteration is linear, we can write down a linear mapping that maps
${\bf u}_1^k$ onto ${\bf u}_1^{k+1}$ for arbitrary $k$. Let this have
Lipschitz constant $L_1$, then we have
\[\|\epsilon_k\|\leq \tau_r\kappa({\bf A}) L_1^k \|{\bf u}_1^0-{\bf u}_1^1\|.\]
Thus the perturbation has limit zero if $L_1<1$. This is the case if
and only if the sequence $({\bf u}_1^k)_k$ is convergent, which is in
fact the case provided that $f(x,y)$ is sufficiently harmless, as can
be seen from the literature on domain decomposition methods,
e.g. \cite[ch. 4]{quaval:99}. 

Analagously to the Picard iteration, if we choose the absolute termination criterion \eqref{linabsterm}
or the relative one based on the right hand side \eqref{linreltermb}, we obtain a
bound of the form
\[\|\epsilon_k\|\leq \|{\bf A}^{-1}\|\tau_r\|{\bf b}({\bf u}_{\Gamma}^k)\|,\]
respectively
\[\|\epsilon_k\|\leq \|{\bf A}^{-1}\| \tau_a.\]
Again, we cannot make a statement on the limit of $\epsilon_k$. 

In the second case, meaning the iteration with Neumann data
\eqref{deltabound}, we obtain 
\[\|\delta_k\| =\|{\bf u}_{2_{\epsilon}}^{k+1}-{\bf B}^{-1}{\bf
  b}_2({\bf u}_{1_{\epsilon}}^{k+1})\| \leq \|{\bf B}^{-1}\|\|{\bf B}{\bf u}_{2_{\epsilon}}^{k+1}-{\bf
  b}_2({\bf u}_{1_{\epsilon}}^{k+1})\|\]
and analogous arguments produce the same
results for $\delta_k$. Thus by theorem \ref{fsitheorem}, we have that when using the relative
criterion we obtain convergence to the exact solution for any
$\tau_r$.

\subsection{A note on the termination criterion and convergence speed}

It is important to note that under the assumptions, all sequences
considered, wether perturbed or not, are convergent and therefore, the fixed point iteration
will terminate when using the standard criterion
\begin{equation}\label{}
\|\x^{k+1}-\x^k\| \leq TOL.
\end{equation}
However, as just shown, the perturbed iteration converges to an
approximation of the unperturbed fixed point and thus, the algorithm
can terminate when we are in fact not $TOL$-close to the solution. 

A further difference between the different iterations that should be
stressed is that the iterations perturbed by a constant are fixed
point iterations, wheras the schemes with a variable perturbation are
in fact, not. Thus, the convergence speed, which is otherwise linear
with constant $L$, is not clear and numerical
evidence suggests that it is in fact slower than for the other
iteration. 

Thus, we could argue to employ the schemes with constant perturbation,
measure the Lipschitz constant numerically after a few iterations and then adjust the tolerance based on
theorem \ref{ortrhetheorem} or \ref{fsitheorem}. Unfortunately, it
is not clear what the $\epsilon$ from these theorems is, respectively,
it is based on quantities that are hard to measure like $\|{\bf
  A}^{-1}\|$. Thus, we
cannot guarantee a certain iteration error in this way, to do this we must employ the nonstandard relative termination criterion. 

Finally, it is important to note that this analysis is mostly relevant
to the time independent case. Otherwise, when considering this inside
an implicit time integration scheme, additional requirements on the
solutions appear, namely that the solutions in the subdomains have a
certain accuracy. 



\section{Numerical Results}

For all numerical experiments, the fixed point iteration is terminated
when the norm $\|\x^{k+1}-\x^k\|_2$ is smaller than
$10^{-14}$. Furthermore, with the exception of the results on the
Picard iteration, all computations were performed in MATLAB, where MATLAB 2012a was used for all computations with the exception of the results in section 3.2.3, where MATLAB 2013a was employed. 

\subsection{Direct Perturbation}



\subsubsection{Testcase: Scalar nonlinear system}

\begin{table}
\centering
\begin{tabular}{c|ccc}
\backslashbox{$L$}{$\epsilon$} & 1e-1 & 1e-2 & 1e-3 \\ \hline
0.009868 & 1.010e-1 & 1.010e-2 & 1.010e-3 \\
0.101239 & 1.090e-1 & 1.089e-2 & 1.089e-3 \\
0.899524 & 2.016e-1 & 1.827e-2 & 1.813e-3 \\
0.996035 & 2.290e-1 & 1.981e-2 & 1.961e-3 \\ \hline
\end{tabular}
\caption{\label{problem2xepsminusxstar}$|x_{\epsilon}-x^*|$ (left) and $\frac{\epsilon}{1-L}$ (right) for different values of
  $\epsilon$ and $L$ for the solution of the scalar nonlinear equation \eqref{scalarnonlineareq}}
\end{table}
As a first example, we employ the nonlinear scalar equation
\begin{equation}\label{scalarnonlineareq}
x=e^{\gamma x}/4
\end{equation}
with $x\in [0,1]$ and $\gamma<1$ given. Thus, the Lipschitz constant $L$ on $[0,1]$ is equal to $\gamma e^{\gamma}/4<1$. We solve this equation for $\gamma=0.3, 1.145, 1.2$. 

Employing the fixed point method with constant perturbation, we
provide the values of $|x_{\epsilon}-x^*|$ in table
\ref{problem2xepsminusxstar}. The difference in solutions is larger than one and proportional to $\epsilon$, as suggested by theorem \ref{ortrhetheorem}. However, the
dependence on the Lipschitz constant is very weak and the error does
not become worse when it approaches one. This is because the
problem is nonlinear and thus, the Lipschitz constant is domain
dependent. The local Lipschitz constant near the solution is actually
well smaller than one, which reminds us that for nonlinear problems, the
Lipschitz constant does not always describe a problem well. 

Furthermore, we tested the adaptive strategy and there $|x_{\epsilon}-x^*|$ tends to machine accuracy, as predicted by the theory. 

\subsubsection{Testcase: Picard iteration}
\label{picardresultssection}
We now consider the Picard iteration \eqref{picardfixed}. The
equation system considered arises from the discretization of the
incompressible Navier-Stokes equations on the unit square with
homogeneous Dirichlet boundary conditions and a viscosity of $\nu
=1/1000$. The grid is cartesian with $128 \times 128 = 16384$
cells. For the computations, the code MooNMD
by John et. al. \cite{johmat:04} was used. The Finite Element discretization employs $Q_21$/$P_1$ elements, resulting in 181250 unknowns overall, thereof 66049 for each velocity component (including the Dirichlet nodes) and 49152 for the pressure.

The right hand side is chosen that the solution is given by
$(u_1,u_2) = (d \psi/d y, - d \psi/ d x)$ with $\psi = x^2 (1-x)^2 y^2 (1-y)^2$, resulting in 

\[u_1(x,y) = x^2 (1-x)^2 [ 2 y (1-y)^2 - 2 y^2 (1-y) ],\]
\[u_2(x,y) = [ 2 x^2 (1-x) - 2 x (1-x)^2 ] y^2 (1-y)^2,\]
\[p(x,y) = x^3 + y^3 - 1/2.\]

As initial guess for the Picard iteration, the zero vector is used. To solve the linear systems \eqref{picardlineq}, GMRES is employed where the initial guess is the current Picard iterate. The Picard iteration is terminated either when the quantity
\[\|res\|=\|{\bf A}(\u^{n+1})\u^{n+1}-{\bf b}\|\]
is smaller than a tolerance $TOL$ or when GMRES terminates immediately without performing an iteration, implying that $\u^{n+1}=\u^n$. 

\begin{table}
\centering
\begin{tabular}{cc|ccc}
TOL & $\tau_r$ & Fixp. iter & GMRES it. & $\|res\|$ \\ \hline
\multirow{7}{*}{1e-14} & 1e-01  & 12 & 42 & 1.793e-15 \\
& 1e-02  & 9   & 48 & 8.559e-15 \\
& 1e-03  & 9   & 58 & 2.318e-15 \\
& 1e-04  & 8   & 57 & 2.121e-15 \\
& 1e-05  & 7   & 59 & 5.299e-15 \\
& 1e-06  & 7   & 64 & 5.266e-15 \\
& 1e-07  & 7   & 67 & 5.265e-15 \\ \hline
1e-07 & 1e-01 & 4 & 11 & 5.66989e-08 \\ \hline
\end{tabular}
\caption{\label{picardproblemreltol}$\|res\|$,
  total number of GMRES iterations and fixpoint iterations for different values of
  $\tau_r$ when using termination criterion \eqref{linrelterm}}
\end{table}
In the first block of table \ref{picardproblemreltol}, $\|res\|$, as well as the total number
of inner GMRES iterations and the number of Picard iterations needed to reach machine accuracy (TOL=$10^{-14}$) are
shown for different values of the relative tolerance $\tau_r$ in GMRES, where the termination criterion \eqref{linrelterm} was used. As predicted by the
theory, all schemes converge to the exact solution. Furthermore, it
takes slightly more fixed point iterations to reach machine accuracy
if the linear systems are solved very inaccurately. Nevertheless, the most
efficient scheme is the one with $\tau_r=1e-1$. 

\begin{table}
\centering
\begin{tabular}{cc|ccc}
TOL & $\tau_a$ & Fixp. iter & GMRES it. & $\|res\|$ \\ \hline
\multirow{8}{*}{1e-14} & 1e-01  & 1 &  0 & 8.212092e-03 \\
& 1e-02  & 1 &  0 & 8.212092e-03 \\
& 1e-03  & 1 &  4 & 6.356705e-04 \\
& 1e-04  & 1 &  6 & 2.131814e-05 \\
& 1e-05  & 1 &  7 & 4.018437e-06 \\
& 1e-06  & 2 & 11 & 3.016563e-07 \\
& 1e-07  & 2 & 14 & 6.249272e-08 \\
& 1e-14  & 7   & 70 & 5.916e-15 \\ \hline
\end{tabular}
\caption{\label{picardproblemabstol}$\|res\|$,
  total number of GMRES iterations and fixpoint iterations for different values of
  $\tau_a$ when using termination criterion \eqref{linabsterm}}
\end{table}
In table \ref{picardproblemabstol}, we show the same quantities, but
for the termination criterion \eqref{linabsterm}. As can be seen, the Picard
iteration does not converge to the exact solution, and how close we
get is proportional to $\tau_a$. This suggests that here, the upper
bound on the perturbation is accurate, thus having a situation as in
theorem \ref{ortrhetheorem}. 

To illustrate the difference in efficiency, we can compare the last row of table \ref{picardproblemabstol} with the first row of table \ref{picardproblemreltol} which tells us that to reach machine accuracy, the method using the relative criterion is faster. A more realistic is obtained by choosing a less strict tolerance. Thus we performed a computation with $TOL=10^{-7}$ and $\tau_r=0.1$ for the first method, the result of which can be seen in the last line of table \ref{picardproblemreltol}. To obtain the same accuracy with the second method, $\tau_a$ has to be chose as $1e-0.7$. Again, the first method is slightly more accurate than the second.

\subsection{Nested fixed point iteration}

\subsubsection{Testcase: Linear Equation System with Matrix Product}
We now consider the linear problem
\begin{equation}\label{linearfsieq}
({\bf I}-{\bf AB}){\bf x}={\bf b} \Leftrightarrow {\bf x}={\bf ABx}+{\bf b}
\end{equation}
with
\[{\bf A}=\left( \begin{array}{cc} \alpha & 0 \\ 0.001 & 0.001 \end{array}
\right ), \, {\bf B}=\left( \begin{array}{cc} \beta & 0 \\ 0.001 & 0.001 \end{array}
\right ), \, {\bf b}=\left ( \begin{array}{c} 1\\ 1 \end{array} \right
).\]
Thus, $\S(\x)={\bf Ax}+{\bf b}$ with $L_S=\|{\bf A}\|_2\approx \alpha$ and ${\bf F}(\x)={\bf
  Bx}$ with $L_F=\|{\bf B}\|_2\approx \beta$.

\begin{table}[hbt]
\centering
\begin{tabular}{c|l|ccc}
$\epsilon$ & \backslashbox{$\alpha$}{$\beta$}  & 0.1 & 0.9 & 0.99 \\ \hline
\multirow{6}{*}{1e-1} & \multirow{2}{*}{0.1}  & 1.111e-1 & 1.209e-1 & 1.221e-1 \\
                      &                       & 1.058e-1 & 1.111e-1 & 1.117e-1 \\ \hline
                      & \multirow{2}{*}{0.9}  & 2.087e-1 & 1.000e-0 & 1.743e-0 \\
                      &                       & 1.638e-1 & 7.107e-1 & 1.235e-0 \\ \hline
                      & \multirow{2}{*}{0.99} & 2.209e-1 & 1.826e-1 & 1.000e+1 \\ 
                      &                       & 1.715e-1 & 1.293e-1 & 7.071e-0 \\ \hline \hline
\multirow{6}{*}{1e-2} & \multirow{2}{*}{0.1}  & 1.111e-2 & 1.209e-2 & 1.221e-2 \\
                      &                       & 1.058e-2 & 1.111e-2 & 1.117e-2 \\ \hline
                      & \multirow{2}{*}{0.9}  & 2.088e-2 & 1.000e-1 & 1.743e-1 \\
                      &                       & 1.638e-2 & 7.107e-2 & 1.235e-1 \\ \hline
                      & \multirow{2}{*}{0.99} & 2.209e-2 & 1.826e-1 & 1.000e-0 \\
                      &                       & 1.715e-2 & 1.293e-1 & 7.071e-1 \\ \hline\hline
\multirow{6}{*}{1e-3} & \multirow{2}{*}{0.1}  & 1.111e-3 & 1.209e-3 & 1.221e-3 \\
                      &                       & 1.058e-3 & 1.111e-3 & 1.117e-3 \\ \hline
                      & \multirow{2}{*}{0.9}  & 2.088e-3 & 1.000e-2 & 1.743e-2 \\
                      &                       & 1.638e-3 & 7.107e-3 & 1.235e-3 \\ \hline
                      & \multirow{2}{*}{0.99} & 2.209e-3 & 1.826e-2 & 1.000e-1 \\ 
                      &                       & 1.715e-3 & 1.293e-3 & 7.071e-2 \\ \hline
\end{tabular}
\caption{\label{problem3xepsminusxstar}Estimate \eqref{fsiestimateepsdelta} and $\|\x_{\epsilon}-\x^*\|_2$ for different values of
  $\epsilon$, $L_S$ and $L_F$ for equation \eqref{linearfsieq}}
\end{table}
As a perturbation, we use a constant vector with eucledian norm $\epsilon=\delta$. The difference $\|\x_{\epsilon}-\x^*\|_2$ can be seen in table
\ref{problem3xepsminusxstar}. As initial guess, the zero vector was
used. The results demonstrate that \eqref{fsiestimateepsdelta} is a very good estimate of the true error and that the errors are perfectly proportional to $\epsilon$. We furthermore tested the adaptive strategy and that iteration indeed converges to $\x^*$.

\subsubsection{Testcase: Scalar nonlinear system}
As a second example, we employ the nonlinear scalar problem
\begin{equation}\label{scalarnonlinearfsieq}
x=0.25 \gamma_1 e^{\gamma_2 x^2}
\end{equation}
with $x\in [0,1]$, $S(x)=0.25 \gamma_1 e^x$ and $L_S=0.25 \gamma_1 e$ and $F(x)=\gamma_2 x^2$ and $L_F=2 \gamma_2$. 

\begin{table}
\centering
\begin{tabular}{c|l|cccc}
$\epsilon$ & \backslashbox{$L_S$}{$L_F$} & 0.01 & 0.1 & 0.9 & 0.99 \\ \hline
\multirow{9}{*}{1e-1} 
                       & \multirow{3}{*}{0.1}  & 1.101e-1 & 1.111e-1 & 1.209e-1 & 1.221e-1 \\
                       &                       & 1.038e-1 & 1.038e-1 & 1.039e-1 & 1.040e-1 \\
                       &                       & 1.039e-1 & 1.039e-1 & 1.042e-1 & 1.042e-1 \\ \hline
                       & \multirow{3}{*}{0.9}  & 1.917e-1 & 2.088e-1 & 1.000e-0 & 1.743e-0 \\
                       &                       & 1.463e-1 & 1.485e-1 & 1.724e-1 & 1.760e-1 \\
                       &                       & 1.350e-1 & 1.370e-1 & 1.618e-1 & 1.658e-1 \\ \hline
                       & \multirow{3}{*}{0.99} & 2.010e-1 & 2.209e-1 & 1.826e-0 & 1.000e+1 \\ 
                       &                       & 1.527e-1 & 1.555e-1 & 1.896e-1 & 1.949e-1 \\
                       &                       & 1.386e-1 & 1.411e-1 & 1.746e-1 & 1.806e-1 \\ \hline \hline
\multirow{9}{*}{1e-2} 
                       & \multirow{3}{*}{0.1}  & 1.101e-2 & 1.111e-2 & 1.209e-2 & 1.221e-2 \\
                       &                       & 1.038e-2 & 1.038e-2 & 1.040e-2 & 1.040e-2 \\
                       &                       & 1.037e-2 & 1.037e-2 & 1.038e-2 & 1.039e-2 \\ \hline
                       & \multirow{3}{*}{0.9}  & 1.917e-2 & 2.088e-2 & 1.000e-1 & 1.743e-1 \\
                       &                       & 1.463e-2 & 1.485e-2 & 1.725e-2 & 1.760e-2 \\
                       &                       & 1.334e-2 & 1.350e-2 & 1.525e-2 & 1.551e-2 \\ \hline
                       & \multirow{3}{*}{0.99} & 2.010e-2 & 2.209e-2 & 1.826e-1 & 1.000e-0 \\
                       &                       & 1.527e-2 & 1.556e-2 & 1.896e-2 & 1.949e-2 \\
                       &                       & 1.368e-2 & 1.388e-2 & 1.621e-2 & 1.658e-2 \\ \hline\hline
\multirow{9}{*}{1e-3} 
                       & \multirow{3}{*}{0.1}  & 1.101e-3 & 1.111e-3 & 1.209e-3 & 1.221e-3 \\
                       &                       & 1.038e-3 & 1.038e-3 & 1.039e-3 & 1.040e-3 \\
                       &                       & 1.037e-3 & 1.037e-3 & 1.038e-3 & 1.038e-3 \\ \hline
                       & \multirow{3}{*}{0.9}  & 1.917e-3 & 2.088e-3 & 1.000e-2 & 1.743e-2 \\
                       &                       & 1.463e-3 & 1.485e-3 & 1.725e-3 & 1.760e-3 \\
                       &                       & 1.333e-3 & 1.348e-3 & 1.518e-3 & 1.542e-3 \\ \hline
                       & \multirow{3}{*}{0.99} & 2.010e-3 & 2.209e-3 & 1.826e-2 & 1.000e-1 \\ 
                       &                       & 1.527e-3 & 1.555e-3 & 1.896e-3 & 1.949e-3 \\
                       &                       & 1.366e-3 & 1.386e-3 & 1.612e-3 & 1.646e-3 \\ \hline
\end{tabular}
\caption{\label{problem4xepsminusxstar}Global estimate, local estimate and $|x_{\epsilon}-x^*|$ for different values of
  $\epsilon$, $L_S$ and $L_F$ for equation \eqref{scalarnonlinearfsieq}}
\end{table}
The initial guess in the following numerical experiments is $x^0=0.5$. In table \ref{problem4xepsminusxstar}, we show several quantities for different values of $\epsilon$, $L_S$ and $L_F$ where a constant perturbation $\epsilon=\delta$ is employed. First, the estimate \eqref{fsiestimateepseqdelta} using the Lipschitz constants on the interval [0,1], which is referred to as the global estimate. Then the local estimate, which is \eqref{fsiestimateepseqdelta} using the derivatives in the solution, giving an estimate of a local Lipschitz constant. This is reasonable, since all functions are monotonic. Finally, the difference $|x_{\epsilon}-x^*|$ itself. 

As can be seen, we again have the proportionality to
$\epsilon$. Furthermore, we see that only when both $L_F$ and $L_S$
are close to one, an influence on the error can be observed, as
suggested by theorem \ref{fsitheorem}. Finally, we test the adaptive
strategy as an example of perturbations converging to zero and again, we obtain convergence of the new sequence to $x^*$.

\subsubsection{Testcase: Transmission Problem}

We now consider the transmission problem \eqref{transmissioneq}. Specifically, we use $\Omega_1=[0,1]\times [0,1]$, $\Omega_2=[1,2]\times [0,1]$ and
\begin{eqnarray*}
f(x,y)=&\sin \pi y^2(\pi \cos \frac{\pi}{2}x^2-\pi^2 x^2 \sin
\frac{\pi}{2}x^2) \\
&+ \sin\frac{\pi}{2}x^2(2\pi \cos \pi y^2 -4\pi^2y^2
\sin\pi y^2).
\end{eqnarray*}
This was chosen such that the solution is 
\begin{equation}
u(x,y)=\sin \pi y^2 \sin \frac{\pi}{2} x^2,
\end{equation}
which satisfies the boundary conditions. 

\begin{figure}[ht]
\includegraphics[width=7cm]{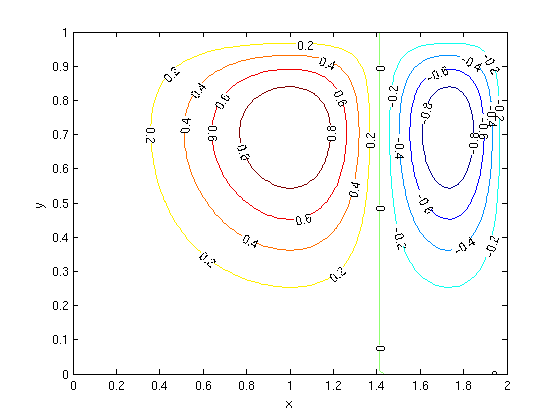}
\includegraphics[width=7cm]{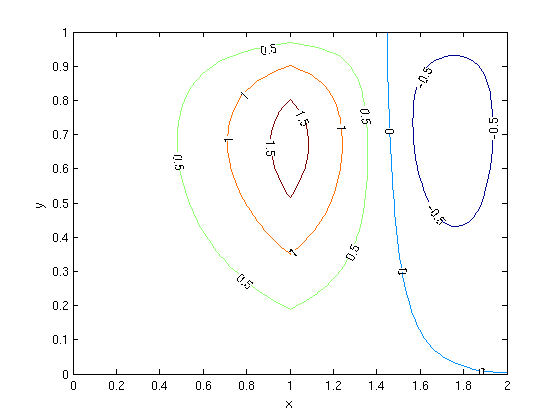}
\caption{\label{transmissionsolutionfigure}Exact and discrete solution with $\Delta x=1/40$}
\end{figure}
We discretize this problem using central differences with a constant mesh width
of $\Delta x=\Delta y$. As initial guess for
the Dirichlet-Neumann procedure, we employ a vector of all zeros. All
linear systems are solved using CG. The exact solution and the discrete solution with $\Delta x=1/40$ can be
seen in figure \ref{transmissionsolutionfigure}. 

\begin{table}[hbt]
\centering
\begin{tabular}{c|cccc}
\backslashbox{$\tau$}{$\Delta x$} & 1/10 & 1/20 & 1/40 & 1/80 \\ \hline
1e-1 & 7.606e-1 & 4.189e-0 & 3.440e+149 & 3.755e+148 \\
1e-2 & 9.620e-2 & 2.502e-2 & 2.621e-1 & 7.288e-1 \\
1e-3 & 1.230e-2 & 2.773e-3 & 1.192e-1 & 1.254e-1 \\ 
1e-4 & 9.110e-4 & 1.033e-3 & 1.074e-2 & 2.602e-2 \\ \hline
\end{tabular}
\caption{\label{transmissionxepsminusxstarrelb}$\|x_{\epsilon}-x^*\|_2$ for different values of
  $\tau$ and $\Delta x$ for the transmission problem
  \eqref{transmissioneq} with relative termination criterion \eqref{linreltermb}}
\end{table}

\begin{table}[hbt]
\centering
\begin{tabular}{c|cccc}
\backslashbox{$\tau$}{$\Delta x$} & 1/10 & 1/20 & 1/40 & 1/80 \\ \hline
1e-1 & 6.643e-3 & 7.308e-3 & 8.993e-3 & 7.187e-3 \\
1e-2 & 7.727e-4 & 6.344e-4 & 8.048e-4 & 6.775e-4 \\
1e-3 & 7.117e-5 & 8.603e-5 & 8.354e-5 & 6.880e-5 \\ 
1e-4 & 5.497e-6 & 7.426e-6 & 7.685e-6 & 6.470e-6 \\ \hline
\end{tabular}
\caption{\label{transmissionxepsminusxstarabs}$\|x_{\epsilon}-x^*\|_2$ for different values of
  $\tau$ and $\Delta x$ for the transmission problem
  \eqref{transmissioneq} with absolute termination criterion \eqref{linabsterm}}
\end{table}
We first look at the convergence properties of the fixed point schemes
for different mesh widths and different termination criteria. The
difference $\|\x_{\epsilon}-\x^*\|_2$ for a constant tolerance
$\tau$ in both CG-subsolvers can be seen in table \ref{transmissionxepsminusxstarrelb}
for the relative termination criterion \eqref{linreltermb} and in
table \ref{transmissionxepsminusxstarabs} for the absolute termination
criterion \eqref{linabsterm}. As can be seen, the schemes behave
essentially as if the perturbation were constant and do not converge
to the exact solution. In particular, for the relative termination criterion \eqref{linreltermb}, the error becomes large for smaller mesh widths and is up to a 100 times larger than the desired tolerance. Not that this criterion is the one implemented in the MATLAB version of CG and that thus, a native implementation of the Dirichlet-Neumann iteration in MATLAB will not produce a correct solution. 

Otherwise, there is again a proportionality to $\tau$, though it's not as clear this time. We
attribute this to the fact that a relative tolerance in CG is only an
upper bound on the perturbation, which can in fact be much smaller
than $\tau$ if CG oversolves. Furthermore, the perturbed solutions
become in general less accurate when the mesh is refined. This can be
explained by the dependence on the norms of ${\bf A}^{-1}$ and ${\bf B}^{-1}$, which increase with decreasing mesh width. 

In the case of the termination criterion \eqref{linrelterm}, we
recover the exact solution, as predicted by the theory. 

\begin{table}
\centering
\begin{tabular}{|c|cc|cc|}
\backslashbox{$\tau_r$}{$\Delta x$} & \multicolumn{2}{c|}{1/10} & \multicolumn{2}{c|}{1/20} \\ \hline
     & \#FP & \#CG &  \#FP & \#CG \\
1e-1 & 106 & 2220 & 205 & 8298  \\
1e-2 & 105 & 2903 & 205 & 11537 \\
1e-3 & 105 & 3121 & 205 & 12585 \\
1e-4 & 105 & 3369 & 208 & 13765 \\ \hline
\backslashbox{$\tau_r$}{$\Delta x$} & \multicolumn{2}{c|}{1/40} & \multicolumn{2}{c|}{1/80} \\ \hline
     & \#FP & \#CG & \#FP & \#CG \\
1e-1 & 401 & 29224 & 379 & 40556 \\
1e-2 & 401 & 44321 & 803 & 156774 \\
1e-3 & 399 & 48341 & 759 & 181789 \\
1e-4 & 402 & 53478 & 835 & 222359 \\ \hline
\end{tabular}
\caption{\label{fsiproblemiterations}Total CG iterations for the transmission problem \eqref{transmissioneq} for different numbers of $\tau_r$ and $\Delta x$ when using termination criterion \eqref{linrelterm} and $TOL=10e-14$.}
\end{table}
We now consider the total number of CG and fixed point iterations when using  the
termination criterion \eqref{linrelterm} for different tolerances
$\tau$ and different mesh widths $\Delta x$. As can be seen in table
\ref{fsiproblemiterations}, the number of CG iterations increases with
decreasing mesh width, which is well known behavior due to the
spectrum getting more widely distributed on the real
line. Furthermore, the number of fixed point iterations is almost
independent of $\tau_r$. Thus, the most efficient variant is to solve
the systems only up to $\tau_r=1e-1$. 

Finally, we compare the different termination criteria for values of TOL more relevant in practice. Hereby, we assume that the user wants to have a solution that is TOL close to the exact one. Based on the theory discussed here, there are three choices: Using \eqref{linrelterm} with $\tau_r=1e-1$ independent of TOL, using \eqref{linabsterm} with $\tau_a=TOL$ or \eqref{linreltermb} with $\tau_r=TOL$, meaning that for the latter ones, we have to solve the inner iteration more accurately the more accurate we want the outer one. Hereby, we choose the largest problem with $\Delta x=1/80$. 
\begin{table}
\centering
\begin{tabular}{|c|cccc|cccc|}
 & \multicolumn{4}{c|}{Relative crit. \eqref{linrelterm}} & \multicolumn{4}{c|}{Relative crit. \eqref{linreltermb}} \\ \hline
TOL  & $\tau_r$ & $$\#FP & \#CG & $\|{\bf x}-{\bf x}^*\|$ & $\tau_r$ & \#FP & \#CG & $\|{\bf x}-{\bf x}^*\|$ \\ \hline
1e-1 & 1e-1 & 22 & 2905 & 8.215e-1   & 1e-1 & - & - & - \\
1e-2 & 1e-1 & 43 & 5258 & 1.203e-1   & 1e-2 & 95 & 10666 & 3.225e-0\\
1e-3 & 1e-1 & 70 & 8256 & 1.397e-2   & 1e-3 & 64 & 10227 & 5.588e-1 \\
1e-4 & 1e-1 & 106 & 12204 & 9.919e-4 & 1e-4 & 116 & 22769 & 2.602e-2 \\ \hline
 & \multicolumn{4}{c|}{Absolute crit. \eqref{linabsterm}} \\ \hline
TOL  & $\tau_a$ & $$\#FP & \#CG & $\|{\bf x}-{\bf x}^*\|$ \\ \hline
1e-1 & 1e-1 & 35 & 11459 & 2.476e-0   \\
1e-2 & 1e-2 & 94 & 30550 & 2.454e-1   \\
1e-3 & 1e-3 & 153 & 50242 & 2.441e-2   \\
1e-4 & 1e-4 & 212 & 70848 & 2.428e-3 \\ \hline
\end{tabular}
\caption{\label{fsiproblemefficiency}Total CG iterations for the transmission problem \eqref{transmissioneq} for different numbers of $TOL$, $\tau_r$ and $\tau_a$ when using different termination critera.}
\end{table}

The results are depicted in table \ref{fsiproblemefficiency}. For the computation with a - divergence was observed. Otherwise, the schemes roughly obey the desired behavior that the error is proportional to TOL, although all are above the desired accuracy. Otherwise, the scheme corresponding to \eqref{linrelterm} is about a factor of five faster than that corresponding to \eqref{linabsterm} and up to a factor of two faster than that corresponding to \eqref{linabsterm} while providing more accurate results. 

\section{Summary and Conclusions}

We considered perturbed fixed point iterations where the perturbation
results from inexact solves of equation systems by iterative
solvers. Thereby, we extended a perturbation result for fixed point
equations to the case of a nested fixed point equation. Applying these
results to the Picard- and the Dirichlet-Neumann iteration for steady states, we showed that these
converge to the exact solution indepently of the tolerance in the subsolver, if a specific
relative termination criterion is employed. This justifies extremely
coarse solves in the inner solvers and suggests the use of GMRES as Krylov subspace solver for unsymmetric systems. 

If an absolute
or standard relative criterion is used, the theory indicates that we
will not converge to the exact solution. Numerical results demonstrate
this behavior. 

Thus, to obtain a certain accuracy in the fixed point solution when using a standard relative or absolute criterion, we have to solve the inner systems more accurate the tighter the tolerance, whereas for the nonstandard relative criterion we can solve the inner systems very coarsely independent of desired accuracy. Numerical results show that this is the most efficient way to treat these systems.

\section*{Acknowledgements}

Part of this work was funded by the German Research Foundation (DFG)
as part of the collaborative research area SFB TRR 30, project
C2. Furthermore, I'd like to thank Gunar Matthies for performing the numerical
experiments on the Picard iteration in section 3.1.2.


\end{document}